\input amsppt.sty

\def\phi{\varphi}
\def\epsilon{\varepsilon}

\NoBlackBoxes
\magnification=1200

\topmatter
\title
On Bergman completeness of non-hyperconvex domains
\endtitle
\author
Marek Jarnicki, Peter Pflug \& W\l odzimierz Zwonek
\endauthor
\address
Instytut Matematyki, Reymonta 4, 30-059 Krak\'ow, Poland
(first and third author)
\endaddress
\email
jarnicki\@im.uj.edu.pl,\;
zwonek\@im.uj.edu.pl
\endemail
\address
Carl von Ossietzky Universit\"at Oldenburg,
Fachbereich 6 -- Mathematik, Postfach 2603, 26111 Oldenburg, Germany
(second author)
\endaddress
\email
pflug\@mathematik.uni-oldenburg.de\endemail
\thanks
The paper has been supported by the KBN grant
No 2 P03A 017 14 and by Nieders\"achsisches Ministerium f\"ur
Wissenschaft und Kunst.
\endthanks
\abstract
In the paper we study the problems of the boundary behaviour
of the Bergman kernel and the Bergman completeness in some
classes of bounded pseudoconvex domains, which contain
also non-hyperconvex domains. Among the classes for which
we prove the Bergman completeness
and the convergence of the Bergman kernel
to infinity while tending to the boundary
are all bounded pseudonvex balanced
domains, all bounded Hartogs domains with balanced fibers
over regular domains and some bounded Laurent-Hartogs domains.
\endabstract
\endtopmatter

\document
\subheading{0. Introduction} The aim of the paper is to present some new
results concerning Bergman completeness and the boundary behaviour
of the Bergman kernel in bounded pseudoconvex but not necessarily hyperconvex
domains. We are interested
in the following exhausting property of the Bergman kernel:
$$
K_D(z)\to\infty\;\text{ as $z\to\partial D$}.\tag{$\ast$}
$$
The starting point for our considerations
may be the following two recent results:

-- any bounded hyperconvex domain
satisfies \thetag{$\ast$} (see \cite{Ohs~2}),

-- any bounded hyperconvex domain is Bergman complete
(see \cite{Blo-Pfl} and \cite{Her}).

Both properties mentioned above are closely related. In particular,
the Bergman completeness is often proved after proving the property
\thetag{$\ast$}. To the best of our knowledge there are no known
examples of bounded Bergman complete domains not satisfying
\thetag{$\ast$}.

The existence of non-hyperconvex bounded domains with \thetag{$\ast$}
is very well-known and easy (take the Hartogs triangle).
On the other hand the existence of
bounded pseudoconvex but non-hyperconvex
Bergman complete domains is not so trivial but also known (see \cite{Chen},
\cite{Her}, \cite{Zwo}).

In our paper we shall present a class of domains satisfying
the above properties. The classes of domains which we consider
are the following:  bounded pseudoconvex balanced domains,
Hartogs domains with $m$-dimensional balanced fibers,
Hartogs-Laurent domains and Zalcman
type domains (domains in the unit disc
with complements composing of infinitely many closed discs).

Among others we prove the following results.

All bounded pseudoconvex balanced domains satisfy \thetag{$\ast$}
and are Bergman complete.
The latter result gives the positive answer
to the question posed in \cite{Jar-Pfl~1} and \cite{Jar-Pfl~2}.
Note that if the Minkowski functional of the considered domain is continuous
then the domain is hyperconvex and the result follows
from the above mentioned theorems.

Any bounded pseudoconvex Hartogs domain with $m$-dimensional balanced fibers
over a domain with the property
\thetag{$\ast$} satisfies \thetag{$\ast$}. Any bounded
pseudoconvex Hartogs domain over a $c^{i}$-complete domain
(which implies automatically Bergman completeness) is Bergman complete.
In particular, there are bounded and pseudoconvex non-fat domains
that are Bergman complete and satisfy \thetag{$\ast$}.

On the other hand we show that there are bounded fat domains in $\Bbb C$
(some Zalcman type domains)
not satisfying \thetag{$\ast$} -- this gives an answer to a question posed
in \cite{Jar-Pfl~2}.

\subheading{1. Definitions and known results}
Let us denote by $E$ the unit disc in $\Bbb C$.

Let $D$ be a
bounded domain in $\Bbb C^n$. Let us denote by
$L_h^2(D)$ square integrable holomorphic functions on $D$. $L_h^2(D)$
is a Hilbert space with the scalar product induced from $L^2(D)$.
Let us define the {\it Bergman kernel of $D$}
$$
K_D(z)=\sup\{\frac{|f(z)|^2}{||f||_{L^2(D)}^2}:f\in L^2_h(D),
f\not\equiv 0\}.
$$
Among other well-known properties let us recall only two of them
(see e.g \cite{Jar-Pfl~2}).

If $D_1\subset D_2$ are bounded domains in $\Bbb C^n$ then
$K_{D_2}(z)\leq K_{D_1}(z)$, $z\in D_1$.

If $\{D_j\}_{j=1}^{\infty}$ is an increasing sequence of domains in $\Bbb C^n$
whose union is a bounded domain $D$,
then $K_{D_j}$ tends decreasingly and locally uniformly to $K_D$.

It is well-known that $\log K_D$ is a smooth plurisubharmonic function.
Therefore, we may define
$$
\beta_D(z;X):=\left(\sum_{j,k=1}^n\frac{\partial^2\log K_D(z)}
{\partial z_j\partial \bar z_k}X_j\bar X_k\right)^{1/2},\;
z\in D,\;X\in\Bbb C^n.
$$
$\beta_D$ is a pseudometric called {\it the Bergman pseudometric}.

For $w,z\in D$ we put
$$
b_D(w,z):=
\inf\{L_{\beta_D}(\alpha)\},
$$
where the infimum is taken over piecewise $C^1$-curves
$\alpha:[0,1]\mapsto D$ joining $w$ and $z$ and
$L_{\beta_D}(\alpha):=\int_0^1\beta_D(\alpha(t);\alpha^{\prime}(t))dt$.

We call $b_D$ {\it the Bergman distance of $D$}.

The Bergman distance (as well as the Bergman metric) is
invariant with respect to biholomorphic mappings.
In other words, for any biholomorphic mapping $F:D\mapsto G$
($D,G\subset\subset \Bbb C^n$) we have
$$
b_G(F(w),F(z))=b_D(w,z),\quad\beta_G(F(w);F^{\prime}(w)X)=\beta_D(w;X),
\;w,z\in D,X\in\Bbb C^n.
$$
A bounded domain $D$ is called
{\it Bergman complete } if any $b_D$-Cauchy sequence is convergent to
some point in $D$ with respect to the standard topology of $D$.

Any bounded Bergman complete domain is pseudoconvex
(see \cite{Bre}).
Let us recall that a bounded domain $D$ is called {\it hyperconvex }
if it admits a continuous negative plurisubharmonic exhaustion function.
Now we may formulate the following very general result:

\proclaim{Theorem 1.1 {\rm (see \cite{B\l o-Pfl}, \cite{Her}, \cite{Ohs~2})}}
Let $D$ be a bounded hyperconvex domain in $\Bbb C^n$.
Then $D$ satisfies \thetag{$\ast$} and $D$ is Bergman complete.
\endproclaim
Our aim is to study the boundary behaviour of the Bergman kernel and the
problem of Bergman completeness. We shall make use of following
powerful tools;
namely, the extension theorem of $L_h^2$-functions, localization principle
of the Bergman kernel and the Bergman metric and criteria for a domain
to be Bergman complete and for the Bergman kernel to tend to infinity
near the fixed point from the boundary.
Let us recall below these results.

\proclaim{Theorem 1.2 {\rm (see \cite{Ohs-Tak})}} Let $D$ be a bounded
pseudoconvex domain in $\Bbb C^n$.
Let $H$ be any affine subspace of $\Bbb C^n$. Then there
is a constant $C\in\Bbb R$ dependent only on diameter of $D$ such
that for any $f\in L_h^2(D\cap H)$ there is $F\in L_h^2(D)$ such that
$||F||_{L^2(D)}\leq C||f||_{L^2(D\cap H)}$.
\endproclaim
In particular, we get from Theorem 1.2 that
$$
K_{D\cap H}(z)\leq\tilde CK_D(z),\;
z\in D\cap H,\tag{1.1}
$$
where $\tilde C\in\Bbb R$ is a constant dependent only on diameter of $D$.

\proclaim{Theorem 1.3 {\rm (see \cite{Die-For-Her}, \cite{Ohs~1})}}
Let $D$ be a bounded pseudoconvex domain in $\Bbb C^n$, $z^0\in\partial D$.
Then for any neighbourhoods $U_1=U_1(z^0)\subset\subset U_2(z^0)=U_2$
there is
a positive constant $C$ such that for any connected component $V$
of $D\cap U_2$ and for any
$z\in U_1\cap V$, $X\in\Bbb C^n$ we have:
$$
\gather
\frac{1}{C}K_V(z)\leq K_D(z)\leq K_V(z),\\
\frac{1}{C}\beta_V(z;X)\leq \beta_D(z;X)\leq \beta_V(z;X).
\endgather
$$
\endproclaim

\proclaim{Theorem 1.4 {\rm (see \cite{Kob}, \cite{Pfl~1})}}
Let $D$ be a bounded domain such that
\thetag{$\ast$} is satisfied and $H^{\infty}(D)$ is dense in $L^2_h(D)$. Then
$D$ is Bergman complete.
\endproclaim

\proclaim{Theorem 1.5 {\rm (see \cite{Pfl~1})}} Let $D$ be a bounded
pseudoconvex domain. Let $z^0\in\partial D$ be such that there
are $r\in (0,1]$, $\epsilon\geq 1$
and a sequence $\{z^{\nu}\}_{\nu=1}^{\infty}$ of points
from $\Bbb C^n\setminus\bar D$
tending to $z^0$
such that $B(z,r||z^{\nu}-z^0||^{\epsilon})\cap D=\emptyset$
(so called 'outer cone condition'). Then
$\lim_{z\to z^0}K_D(z)=\infty$.
\endproclaim

\subheading{2. Balanced domains} Recall that a set $D$
is {\it balanced } if $z\in D$ and $\lambda\in \bar E$ implies
$\lambda z\in D$.

In this section we deal with bounded pseudoconvex balanced domains.
We prove the following result.
\proclaim{Theorem 2.1} Let $D$ be a bounded pseudoconvex balanced domain.
Then $D$ satisfies \thetag{$\ast$} and $D$ is Bergman complete.
\endproclaim
Note that if the Minkowski functional of $D$ is continuous then $D$ is
hyperconvex and the result follows from Theorem 1.1.
Additionally, in this case the theorem has already
been known for a long time (see \cite{Jar-Pfl~1}). Using only a little
more refined methods
than the ones used in the last paper we prove the theorem in general case.
Let us mention here that the problem
whether bounded pseudoconvex balanced domains are Bergman complete
has been stated in \cite{Jar-Pfl~1} and \cite{Jar-Pfl~2}.

\demo{Proof of Theorem 2.1} First we prove the property \thetag{$\ast$}.
Take any point $z^0\in\partial D$.
Fix any $M\in\Bbb R$.
In view of \thetag{1.1} (applied to $H=\Bbb Cz^0$ --
remember that $\Bbb Cz^0\cap D$ is a disc)
there is some $z^1=sz^0$, $0<s<1$
such that $z^1\in D$ and $K_D(z^1)>M$.
It follows from continuity of $K_D$ that
there is some open neighbourhood $U\subset D$ of $z^1$ such that
$K_D(z)>M$ for $z\in U$. Note that for any $z\in U$
the function
$$
u_z:\frac{1}{h(z)}E\owns\lambda\mapsto K_D(\lambda z)
$$
is subharmonic and radial. Therefore, $u_z(t)$,
$0\leq t<\frac{1}{h(z)}$
is increasing (see e.g. \cite{Jak-Jar}).
Consequently, $K_D(z)>M$ for any $z\in([1,\infty)U)\cap D$.
Since $[1,\infty)U$
is a neighbourhood of $z^0$, we finish the proof.

To finish the proof it is sufficient to show that
$H^{\infty}(D)$ is dense in $L^2_h(D)$ (and then use Theorem 1.4).

It is well-known that any holomorphic function $F$
on $D$ is a local uniform limit of a series
$\sum_{k=0}^{\infty}Q_k(z)$,
where $Q_k$ is a homogeneous polynomial of degree $k$
(see e.g. \cite{Jak-Jar}).
Since all $Q_k$ are orthogonal (in $L^2_h(D)$)
and there is an exhausting family of
compact balanced sets of the domain $D$ (on each of them
the functions $Q_k$ are orthogonal),
the standard approximation process leads to the convergence
of $F_N:=\sum_{k=0}^NQ_k$ to $F$ in $L^2(D)$ norm (under the assumption
that $F\in L^2_h(D)$). Since $D$ is bounded, all $F_N$'s are bounded,
which finishes the proof.
\qed
\enddemo

\subheading{3. Hartogs domains}
In the present section
we consider bounded pseudoconvex
Hartogs domains with
$m$-dimensional balanced fibers.
Let $G_D\subset \Bbb C^{n+m}$
denote a bounded pseudoconvex Hartogs domain over $D\subset\Bbb C^n$
with $m$-dimensional balanced fibers, i.e.
$$
G_D=\{(z,w)\in D\times\Bbb C^m:
H(z,w)<1\},
$$
where $D$ is bounded
and pseudoconvex,
$\log H$ is plurisubharmonic
on $D\times\Bbb C^m$, $H(z,\lambda w)=|\lambda|H(z,w)$,
$(z,w)\in D\times\Bbb C^m$,
$\lambda\in\Bbb C$,
and $G_D$ is bounded
(i.e. $H(z,w)\geq C||w||$ for some $C>0$, $z\in D$, $w\in\Bbb C^m$).

Let $G_D$ be as above.
For any $f\in L^2_h(D)$
we define a function
$F(z,w):=f(z)$, $(z,w)\in G_D$. Since $G_D\subset D\times (R E)^m$ for some
$R>0$, we easily get that $||F||_{L^2_h(G_D)}\leq C_1||f||_{L^2_h(D)}$ for
some $C_1>0$ independent of the choice of $f$ (therefore,
$F\in L^2_h(G_D)$). In particular,
we get that
$$
K_D(z)\leq C_2 K_{G_D}(z,0),\; z\in D.\tag{3.1}
$$

\proclaim{Theorem 3.1} Let $G_D$ be a
bounded pseudoconvex Hartogs domain over $D$ with $m$-dimensional
balanced fibers. Fix a point $(z^0,w^0)\in\partial G_D$. Assume
that one of the following three conditions is satisfied:
$$
\gather
z^0\in D,\tag{i}\\
z^0\in\partial D\text{ and } K_D(z)\to\infty \text{ as $z\to z^0$},\tag{ii}\\
\text{there is some neighbourhood $U$ of $(z^0,w^0)$ such that}
\tag{iii}\\
U\cap G_D\subset
\{(z,w)\in\Bbb C^{n+m}:||w||<||z-z^0||^{\delta}\} \text{ for some
$\delta>0$ (in particular, $w^0=0$)}.
\endgather
$$
Then $K_{G_D}(z,w)\to\infty$ as $(z,w)\to (z^0,w^0)$.

In particular, if $D$ satisfies \thetag{$\ast$},
then $G_D$ satisfies \thetag{$\ast$}.
\endproclaim
\demo{Proof}
Consider the case $z^0\in D$ (that is we consider the case
\thetag{i}). Then $H(z^0,w^0)\geq 1$.
In view of the $L^2_h$-extension theorem for any $M\in\Bbb R$
there is $w^1=tw^0$,
$0<t<1$, such that $(z^0,w^1)\in G_D$ and $K_{G_D}(z^0,w^1)>M$. Continuity of
$K_{G_D}$ gives us the existence of an open neighbourhood
$U:=U_1\times U_2$ of $(z^0,w^1)$ in $G_D$ (with $0\not\in U_2$) such that
$K_{G_D}(z,w)>M$ for $(z,w)\in U$. Similarly,
as earlier considering the function
$$
u_{(z,w)}:\frac{1}{H(z,w)}E\owns\lambda\mapsto K_{G_D}(z,\lambda w)
$$
we get a radial subharmonic function such that $u_{(z,w)}(1)>M$, which
gives us that $K_{G_D}(z,w)>M$
for $(z,w)\in (U_1\times[1,\infty)U_2)\cap G_D$.

Consider now the case \thetag{ii}. Then $z^0\in\partial D$.
It follows from \thetag{3.1}
and \thetag{ii} that for any $M\in\Bbb R$
there is some open neighbourhood $U$ of
$z^0$ such that $K_{G_D}(z,0)>M$, $z\in U\cap D$.

Fix $z\in U\cap D$. Fix additionally
for a while $w$ such that $0<H(z,w)<1$. Then the function
$\frac{1}{H(z,w)}E\owns \lambda\mapsto K_{G_D}(z,\lambda w)$
is larger than $M$ at $0$ and is radial and subharmonic;
therefore, increasing. Consequently,
$K_{G_D}(z,w)>M$ for any $w$ with $H(z,w)<1$. Since $z\in U$ was chosen
arbitrarily we get that
$K_{G_D}(z,w)>M$
for any $(z,w)\in G_D$ with $z\in U$.

We are left with the case \thetag{iii}.
Without loss of generality we may assume that $z^0=0$.
Consider points $(0,w_0)\not\in\bar D$ (i.e. $w_0\neq 0$).
Let us consider the balls $B((0,w_0),r||w_0||^{\epsilon})$, where
$\epsilon>0$, $0<r<1$ will be chosen later (independently of $w_0$).
Our aim is to verify that the outer cone condition form Theorem 1.5
is satisfied for a suitable $0<r\leq 1$ and $\epsilon\geq 1$.

Fix $r=\frac{1}{2}$. Consider only $||w_0||<\frac{1}{2}$.
Take a point $(z,w)\in B((0,w_0),\frac{1}{2}||w_0||^{\epsilon})\cap \bar G_D$.
Then $||z||<\frac{1}{2}||w_0||^{\epsilon}$
and $||w_0||-||w||\leq ||w-w_0||\leq
\frac{1}{2}||w_0||^{\epsilon}$.
Consequently, $||w_0||-\frac{1}{2}||w_0||^{\epsilon}\leq
||w||\leq||z||^{\delta}<(\frac{1}{2})^{\delta}||w_0||^{\delta\epsilon}$.
So assuming that $\epsilon$ is large enough
($\epsilon-1>0$, $\epsilon\delta-1>0$, $\delta+\epsilon\delta-2>0$)
we get:
$$
\frac{1}{2}<1-\frac{1}{2}||w_0||^{\epsilon-1}<(\frac{1}{2})^{\delta}
||w_0||^{\epsilon\delta-1}<\frac{1}{2}
$$
-- contradiction. Therefore, in view of Theorem 1.5 we finish the proof.
\qed
\enddemo
The idea of the condition \thetag{iii} comes
from generalizing the phenomenon,
which appears in the Hartogs triangle and the point $(z^0,w^0)=(0,0)$.

It turns out that there are bounded pseudoconvex Hartogs domains
and points from the boundary,
which do not satisfy any from the conditions \thetag{i}--\thetag{iii}
but such that the limes as in Theorem 3.1 exists.

\subheading{Example 3.2} Let $\{a_j\}_{j=1}^{\infty}\subset(0,1)$
be a sequence tending to $0$.
Let us define
$u_k(\lambda)=\log\left(\sum_{j=1}^k(\frac{a_j}{2|\lambda-a_j|})^{n_j}
\right)$,
$\lambda\in E\setminus\{a_1,\ldots,a_k\}$, where $n_j\geq j$.
Note that $u_k(0)<0$, $k=1,2,\ldots$. Define $u:=\lim_{k\to\infty}u_k
=\log\left(\sum_{j=1}^{\infty}(\frac{a_j}{2|\lambda-a_j|})^{n_j}
\right)$
on $E_{\infty}:=E\setminus(\{a_j\}_{j=1}^{\infty}\cup\{0\})$.
The construction ensures us that the sequence
$\{u_k\}_{k=1}^{\infty}$ is locally bounded from above and globally
bounded from below increasing sequence
on $E_{\infty}$ and, therefore,
$u$ is a subharmonic function
on $E_{\infty}$ bounded from below.
Moreover,  $\lim_{x<0, x\to 0}u(x)\leq 0$.
Define $G_{E_{\infty}}:=\{(z,w)\in E_{\infty}\times\Bbb C:\;
|w|<\exp(-u(z))\}$.
Then $G_{E_{\infty}}$ is a bounded pseudoconvex Hartogs domain
with one-dimensional fibers. Note that the point $(0,0)$ does not satisfy
any from the conditions \thetag{i}--\thetag{iii} but one may easily
verify that, choosing if necessary $n_j$ larger, the outer cone condition
from Theorem 1.5 is satisfied (for instance for points $(a_j,a_j)$).
Therefore, the claim of Theorem 3.1 is also satisfied.
Note that $\{(0,0)\}\subsetneq
\partial G_{E_{\infty}}
\cap(\{0\}\times\Bbb C^m)$.

We may prove even more. Namely, the domain $G_{E_{\infty}}$ satisfies
\thetag{$\ast$}. In fact, the points $(z,w)\in\partial G_{E_{\infty}}$,
$z\in\partial E$,
satisfy \thetag{ii}.
The points $(a_k,w)\in\partial G_{E_{\infty}}$ (and then automatically
$w=0$) satisfy \thetag{iii}.
The points
$(z,w)\in\partial G_{E_{\infty}}$,
$z\in E_{\infty}$, satisfy
\thetag{i}. Finally, one may easily verify (proceeding similarily
as in the case $(0,0)$) that
the points $(0,w)\in\partial G_{E_{\infty}}$
satisfy the outer cone condition from Theorem 1.5.

\proclaim{Lemma 3.3} Let $G_D$ be a bounded pseudoconvex Hartogs domain
over $D$ with $m$-dimensional balanced fibers such that
$H^{\infty}(D)$ is dense in $L^2_h(D)$ and,
additionally, assume that there is
some $\epsilon>0$ such that $D\times P(0,\epsilon)\subset G_D$.
Then $H^{\infty}(G_D)$ is dense in $L^2_h(G_D)$.
\endproclaim
\demo{Proof} Take $F\in L^2_h(G_D)$. We know that
$$
F(z,w)=\sum_{\nu=0}^{\infty}F_{\nu}(z,w):=\sum_{\nu=0}^{\infty}
\sum_{\beta\in\Bbb Z_+^m:|\beta|=\nu}f_{\beta}(z)w^{\beta}
$$
where the convergence of $G_N:=\sum_{\nu=0}^NF_{\nu}$ to $F$ is
locally uniform (see e.g. \cite{Jak-Jar}).
Consequently, because of orthogonality
of $w^{\beta}$, similarly as in the proof of Theorem 2.1 the functions
$G_N$ converge in
$L^2(G_D)$ to $F$. It is therefore sufficient to approximate
$f_{\beta}(z)w^{\beta}$ with bounded functions. But because of the assumption
of the lemma one may easily conclude from the Fubini theorem that
$f_{\beta}\in L^2_h(D)$ so $h_N(z)w^{\beta}$,
where $h_N\in H^{\infty}(D)$ and
$h_N\to f_{\beta}$ in $L^2_h(D)$, tends
to $f_{\beta}(z)w^{\beta}$ in $L^2_h(G_D)$.
\qed
\enddemo
\subheading{Remark 3.4} Note that the assumption $D\times P(0,\epsilon)
\subset G_D$ is essential. For instance, $H^{\infty}(E_*)=
H^{\infty}(E)|_{E_*}$
is dense in $L^2_h(E_*)=L^2_h(E)|_{E_*}$
and $H^{\infty}(G_{E_*})$ is not dense in
$L^2_h(G_{E_*})$, where $G_{E_*}$ is the Hartogs triangle,
$G_{E_*}:=\{(z,w)\in E_*\times \Bbb C:|w|<|z|\}$.

\proclaim{Theorem 3.5} Let $G_D$ be a bounded pseudoconvex Hartogs domain
over $D$ with $m$-dimensional balanced fibers. Assume that $D$
satisfies \thetag{$\ast$}, $H^{\infty}(D)$ is dense in $L^2_h(D)$
and there is $\epsilon>0$ such that $D\times P(0,\epsilon)\subset G_D$.
Then $G_D$ is Bergman complete.
\endproclaim
\demo{Proof} Combine Theorem 3.1, Lemma 3.3 and Theorem 1.4.
\qed
\enddemo

Note that Theorem 3.5 cannot be even applied
to arbitrary pseudoconvex bounded Hartogs domain with one dimensional fibers.
However, small change in assumptions on the domain $D$ in Theorem 3.5
will make it possible to prove Bergman completeness of $G_D$
without additional assumptions on the shape of $G_D$.
But before formulating the result we have to introduce the notion
of the inner Carath\'eodory pseudodistance.

For a domain $D\subset\Bbb C^n$ we define {\it the Carath\'eodory-Reiffen
pseudometric}
$$
\gamma_D(z;X):=\sup\{|f^{\prime}(z)X|:\;f\in\Cal O(D,E),\;f(z)=0\},\;
z\in D,\;X\in\Bbb C^n.
$$
{\it The inner Carath\'eodory pseudodistance }
is the integrated form of $\gamma_D$, i.e.
$$
c_D^i(w,z):=\inf\{L_{\gamma_D}(\alpha):\;\alpha:[0,1]\mapsto D
\text{ is a piecewise $C^1$-curve joining $w$ and $z$}\},
$$
where $L_{\gamma_D}(\alpha):=\int_0^1\gamma_D(\alpha(t);
\alpha^{\prime}(t))dt$.
It is well-known that holomorphic mappings are contractions
with respect to $c^i$ (i.e. $c_G^{i}(F(w),F(z))\leq c_D^i(w,z)$
for any $F\in\Cal O(D,G)$, $w,z\in D$).
The last property is not shared by the Bergman
distance (in the class of bounded domains -- see e.g. \cite{Jar-Pfl~2}).
We have additionally that (see e.g. \cite{Jar-Pfl~2})
$$
c_D^i\leq b_D.\tag{3.2}
$$
Exactly as in the case of the Bergman distance we introduce
for bounded domains the notion of $c^i$-completeness.

\proclaim{Theorem 3.6} Let $G_D$ be a bounded pseudoconvex Hartogs domain
over $D$ with $m$-dimensional balanced fibers. Assume that
$D$ is $c^i$-complete.
Then $G_D$ is Bergman complete.
\endproclaim
\demo{Proof} Take any point $(z_0,w_0)\in\partial G_D$. Suppose that
there is a $b_{G_D}$-Cauchy sequence
$\{(z_{\nu},w_{\nu})\}$ tending (in the natural topology of $D$)
to $(z_0,w_0)$.
Because of \thetag{3.2} and the contractivity of $c_D^i$ with respect
to the projection we exclude the case $z^0\in\partial D$.

So assume that $z^0\in D$.
Let $U_1,U_2$ be small open balls with the centre at $z^0$
such that $U_1\subset\subset U_2\subset\subset D$.
There is a sequence of $C^1$-piecewise
curves $\gamma_{\nu,\mu}:[0,1]\mapsto G_D$
such that $\gamma_{\nu,\mu}(0)=(z_{\nu},w_{\nu})$,
$\gamma_{\nu,\mu}(1)=(z_{\mu},w_{\mu})$
and $L_{\beta_{G_D}}(\gamma_{\nu,\mu})<b_{G_D}((z_{\nu},w_{\nu}),
(z_{\mu},w_{\mu}))+\frac{1}{\nu}$, $1\leq\nu<\mu$. We claim that there is
some $\nu_0$ such that $\gamma_{\nu,\mu}([0,1])\subset
G_{U_1}$ ($G_{U_j}:=(U_j\times\Bbb C^m)\cap G_D$, $j=1,2$)
for $\mu>\nu>\nu_0$.
Actually, if it were not the case, then
there would be a sequence of $t_k\in(0,1)$
such that $(u_k,v_k):=\gamma_{\nu_k,\mu_k}(t_k)\not \in G_{U_1}$
(so $u_k\not\in U_1$) and
$b_{G_D}((z_{\nu_k},w_{\nu_k}),(u_k,v_k))\to 0$
as $k$ tends to infinity. But then also
$$
0\leq c_D^i(z_{n_k},u_k)\leq c_{G_D}^i((z_{n_k},w_{n_k}),(u_k,v_k))\to 0
$$
-- contradiction (to the boundedness of $D$).

Note that $G_{U_j}$ satisfies the assumptions of Theorem 3.5,
so $G_{U_j}$ is Bergman complete, $j=1,2$.

Applying the localization principle of the Bergman metric (Theorem 1.3)
we get that
$$
\multline
b_{G_{U_2}}((z_{\nu},w_{\nu}),(z_{\mu},w_{\mu}))\leq
L_{\beta_{G_{U_2}}}(\gamma_{\nu,\mu})\leq\\
C L_{\beta_{G_D}}(\gamma_{\nu,\mu})
\leq C (b_{G_D}((z_{\nu},w_{\nu}),(z_{\mu},w_{\mu}))+\frac{1}{\nu}),\;
\mu>\nu>\nu_0
\endmultline
$$
so $\{(z_{\nu},w_{\nu})\}_{\nu>>1}$ is a $b_{G_{U_2}}$-Cauchy sequence
tending to the boundary of $G_{U_2}$ (in the natural topology of $G_{U_2}$),
which, however,
contradicts the Bergman completeness of $G_{U_2}$.
\qed
\enddemo

\subheading{Remark 3.7} Since any
Kobayashi complete bounded domain
is taut (a bounded domain $D$ in $\Bbb C^n$ is {\it taut } if
for any convergent sequence of mappings $\phi_{\nu}\in\Cal O(E,D)$
its limit $\phi$ satisfies that $\phi(E)\subset D$
or $\phi(E)\subset\partial D$),
there are bounded pseudoconvex balanced domains (in fact
any such that the Minkowski functional is not continuous)
in $\Bbb C^2$ such that no estimate of the type $b_D\leq Ck_D$ holds
(compare \cite{Jar-Pfl~2}).

Note that there are bounded balanced pseudoconvex domains
which are not fat (i.e. $\operatorname{int}(\bar D)\neq D$, see
\cite{Sic}),
so there are Bergman complete domains satisfying \thetag{$\ast$},
which are not fat
(use Theorem 2.1). Other domains having the same property
(but in the class of Hartogs domains) are given below.

Theorem 3.1, Theorem 3.5 and Theorem 3.6
apply among others to the following domain
$$
G_E:=\{(z,w)\in\ E\times\Bbb C:|w|<\exp(-\exp(\sum_{j=1}^{\infty}\alpha_j
\log\frac{|z-a_j|}{2}))\},
$$
where $\alpha_j>0$, $\{a_j\}_{j=1}^{\infty}$ is dense in $E_*$
and $\sum_{j=1}^{\infty}\alpha_j\log|a_j|>-\infty$.
Note that $G_E\subset E^2$, $G_E\neq E^2$
but $\operatorname{int}(\bar G_E)=E^2$.

It follows from Theorem 3.6 that any bounded pseudoconvex Hartogs
domain over a complete bounded pseudoconvex Reinhardt domain
(e.g. the unit disc) is Bergman complete (see \cite{Pfl~2}).

It seems to be natural to ask the question whether Theorem 3.6 remains
true under the assumption that $D$ is Bergman complete.

Since any bounded hyperconvex domain is Bergman complete,
new results concerning Bergman completeness are given in non-hyperconvex case.
In the class of bounded pseudoconvex balanced domains hyperconvexity
is equivalent to tautness and the latter is equivalent to the
continuity of the Minkowski functional associated to the domain.

Below we give a full characterization of tautness and hyperconvexity in the
class of bounded pseudoconvex Hartogs domains
with $m$-dimensional balanced fibers.

\proclaim{Proposition 3.8}
Let $G_D$ be a bounded pseudoconvex Hartogs domain over $D$ with
$m$-dimensional balanced fibers. Then

$G_D$ is taut iff $D$ is taut and $H$ is continuous;

$G_D$ is hyperconvex iff $D$ is hyperconvex and $H$ is continuous.
\endproclaim
\demo{Proof} Note that non-continuity of $H$ gives us the existence
of the sequence $\{(z_{\nu},w_{\nu})\}\subset G_D$
convergent to $(z,w)\in D\times\Bbb C^m$ such that $\lim_{\nu\to\infty}
H(z_{\nu},w_{\nu})=\delta<H(z,w)=1$.
Then the sequence $\phi_{\nu}(\lambda):=
(z^{\nu},\frac{w^{\nu}\lambda}{H(z_{\nu},w_{\nu})})$, $\lambda\in E$,
satisfies $\phi_{\nu}(E)\subset G_D$ and $\phi_{\nu}$ converges
locally uniformly to $\phi$, where
$\phi(\lambda)=(z,\frac{w\lambda}{\delta})$, $\phi(0)\in G_D$
but $\phi(E)\not\subset G_D$, so $G_D$ cannot be taut.

It is trivial to see that tautness (respectively, hyperconvexity) of $G_D$
implies tautness (respectively, hyperconvexity) of $D$.

Hyperconvexity of $D$ delivers us the existence of negative continuous
plurisubharmonic exhaustion function $u$ of $D$. Note that if $H$ is
continuous, then the function $\max\{u(z),\log H(z,w)\}$ is a continuous
negative exhaustion function of $G_D$.

Assume now tautness of $D$ and continuity of $H$.
Consider a sequence $\phi^{\nu}:=(\phi^{\nu}_1,\phi^{\nu}_2)
\in\Cal O(E,G_D)$, which converges locally uniformly to $\phi^0$. Because
of tautness of $D$ either $\phi^0_1\in\Cal O(E,D)$
or $\phi^0_1(E)\subset\partial D$, in the second case
$\phi^0(E)\subset\partial G_D$. So consider the first case.
It easily follows from the maximum principle for subharmonic functions
that either $H(\phi^0(\lambda))=1$ or $H(\phi^0(\lambda))<1$, $\lambda\in E$,
which finishes the proof.
\qed
\enddemo

\subheading{4. Hartogs-Laurent domains}
In this section we consider Hartogs-Laurent domains.
More precisely, let $D$ be a bounded pseudoconvex
domain in $\Bbb C^n$ and let $u,v$
be plurisubharmonic functions on $D$, $u+v<0$ on $D$.
Then we define {\it the Hartogs-Laurent domain  $G:=\{(z,z_{n+1})\in D\times \Bbb C: \exp(u(z))<
|z_{n+1}|<\exp(-v(z))\}$ over $D$}.
We assume additionally that there is some constant $C\in\Bbb R$ such that
$v(z)>C$ (i.e. $G$ is bounded) and $u\not\equiv-\infty$.

\proclaim{Proposition 4.1} Let $G$ be as above (with some $D$, $u$ and $v$).
Assume additionally that $D$ satisfies \thetag{$\ast$}.
Then $G$ satisfies \thetag{$\ast$}.
\endproclaim
\demo{Proof} Since $G\subset\{(z,z_{n+1})\in D\times\Bbb C:
|z_{n+1}|<\exp(-v(z))\}$, we get in view of Theorem 3.1
and because of the contraction property of the Bergman kernel
under inclusion of domains that $K_G(z,z_{n+1})\to\infty$ whenever
$(z,z_{n+1})\to(w,w_{n+1})\in\partial G$, where
$w\in\partial D$ or $|w_{n+1}|\geq\exp(-v(w))$.

Now we consider the case when
$(z,z_{n+1})\to (w,w_{n+1})\in\partial G$, where
$w\in D$ and $|w_{n+1}|\leq\exp(u(w))$.

First we prove that $K_G(z,z_{n+1})\to\infty$ as $(z,z_{n+1})
\to (w,0)\in\partial G$ with $w\in D$. Take a small ball
$U\subset\subset D$ with the centre at $w$. Put $G_U:=G\cap(U\times\Bbb C)$.
We claim that the function $\frac{1}{z_{n+1}}$ is from $L^2_h(G_U)$.
In fact,
$$
\multline
\int_{G_U}\frac{1}{|z_{n+1}|^2}d\Cal L^{2n+2}(z,z_{n+1})=\\
\int_U(\int_{\exp(u(z))<|z_{n+1}|<\exp(-v(z))}\frac{1}{|z_{n+1}|^2}
d\Cal L^2(z_{n+1}))
d\Cal L^{2n}(z)=\\
2\pi\int_U(-v(z)-u(z))d\Cal L^{2n}(z).
\endmultline
$$
Therefore, in view of the local summability of plurisubharmonic functions
(not identical to $-\infty$) the last expression is finite.
Consequently, $K_{G_U}(z,z_{n+1})\to\infty$
as $(z,z_{n+1})\to (w,0)$. And now the localization property
of the Bergman kernel (Theorem 1.3)
implies that
$$
K_G(z,z_{n+1})\to\infty\;\text{ as $(z,z_{n+1})\to(w,0)$.}
$$
We are left with the case $(z,z_{n+1})\to(w,w_{n+1})\in\partial G$,
$w\in D$, $0<\epsilon<|w_{n+1}|\leq\exp(u(w))$. Consider now
the new Laurent-Hartogs domain $G_1$ defined over $D$ with $u$ replaced
by $\tilde u:=\max\{u,\log\epsilon\}$ (and the same $v$). Taking now $\tilde G_1$
to be $\{(z,1/z_{n+1}):(z,z_{n+1})\in G_1\}$ we get that the convergence
of $K_{G_1}(z,z_{n+1})\to\infty$ as $(z,z_{n+1})\to (w,w_{n+1})$
is equivalent to the convergence of $K_{\tilde G_1}(z,z_{n+1})\to\infty$
as $(z,z_{n+1})\to(w,1/w_{n+1})$ (use the invariance of the Bergman kernel
with respect to biholomorphic mappings). Since $\tilde G_1\subset
\{(z,z_{n+1}):z\in D,\;|z_{n+1}|<\exp(-\tilde u(z))\}$
and $1/|w_{n+1}|\geq\exp(-u(w))$
we get (using the contractivity of the Bergman
kernel under inclusion and Theorem 3.1)
that
$K_{G_1}(z,z_{n+1})\to\infty$
as $(z,z_{n+1})\to (w,w_{n+1})$.
And now the localization of the Bergman kernel (Theorem 1.3) implies that
$K_G(z,z_{n+1})\to\infty$ as
$(z,z_{n+1})\to (w,w_{n+1})$.
\qed
\enddemo

\comment
Define for $z\in D$ the following subharmonic functions:
$$
\phi_z:P_z:=\{\exp(u(z))<|\lambda|<\exp(-v(z))\}\owns\lambda
\mapsto K_D(z,\lambda).
$$
Certainly all the functions $\phi_z$ are radial and subharmonic; moreover,
$\phi_z(r)\to\infty$ as $r\to\partial P_z$ -- for $r\to r_0>0$
it is a consequence
of the $L^2_h$ extension theorem (applied to $H:=\{z\}\times\Bbb C$)
and the corresponding convergence in annuli (or the punctured disc)
$P_z$
and for $r\to 0$ it is a consequence of \thetag{4.1}.
Consequently, any from the functions $\phi_z(r)$ is
decreasing in the interval $(\exp(u(z)),r_z)$ and increasing
in the interval $(r_z,\exp(-v(z))$ for some $\exp u(z)<r_z<\exp(-v(z))$.

Now we may complete similarly as in the proof of Proposition 4.1
the proof of the convergence of $K_G(z,z_{n+1})\to\infty$ as
$(z,z_{n+1})\to (w,w_{n+1})$ with $w\in D$, $|w_{n+1}|\leq
\exp(u(w))$.
Fix $M\in\Bbb R$. It follows from continuity of the function $K_G$
and the above mentioned properties of $\phi_z$ that there is some
$R_w<r_w$ and some neighbourhood $U$ of $z$
such that $\phi_z(\lambda)>M$ for any $\exp(u(z))<|\lambda|<R_w$,
$z\in U$, which completes the proof.
\endcomment

\proclaim{Theorem 4.2}  Let $D$ be a bounded
domain pseudoconvex domain
in $\Bbb C^n$, which is $c^i$-complete.
Let $G$ be as above with the additional property that there is some constant $C$
such that $u(z)>C>-\infty$ for any $z\in D$.
Then $G$ is Bergman complete.
\endproclaim
\demo{Proof} We proceed similarly as in the proofs of results in Section 3.
Take any $b_G$-Cauchy sequence $\{(z^{\nu},z_{n+1}^{\nu})\}$ converging
to $(z^0,z_{n+1}^0)\in\partial G$.
One easily excludes (because of $c^i$-completeness of $D$)
the case $z^0\in\partial D$. In case $z^0\in D$ we may exactly as
in the proof of Theorem 3.6 reduce the problem
to the problem of completeness of $\tilde G:=G\cap (U\times\Bbb C)$,
where $U$ is some
small ball around $z^0$, $U\subset\subset D$, such that $U\times A\subset G$,
where $A$ is some annulus. Similarly, as in the proof
of Lemma 3.3, expanding any $L^2_h$-function in the series
$F(z,z_{n+1}):=\sum_{\nu=-\infty}^{\infty}h_{\nu}(z)z_{n+1}^{\nu}$,
$(z,z_{n+1})\in\tilde G$,
we easily get that $F_N(z,z_{n+1}):=\sum_{\nu=-N}^Nh_{\nu}(z)z_{n+1}^{\nu}$
tends in $L^2_h(\tilde G)$ to $F$. Moreover,
$h_{\nu}(z)\in L^2_h(U)$. Since $H^{\infty}(U)$ is dense in $L^2_h(U)$,
we get consequently, that $H^{\infty}(\tilde G)$
is dense in $L^2_h(\tilde G)$ (approximate $h_{\nu}(z)z_{n+1}^{\nu}$
with $h_{\nu,j}(z)z_{n+1}^{\nu}$, where $h_{\nu,j}\in H^{\infty}(U)$
tends to
$h_{\nu}$ in $L^2(U)$, and then use inequality $u>C$),
which in connection with Proposition 4.1
and Theorem 1.4 finishes the proof.
\qed
\enddemo

\subheading{5. Zalcman type domains} In Section 3 we saw that there are
non-fat domains satisfying \thetag{$\ast$} and being Bergman complete. In this
chapter we go into the opposite direction and we find domains, which are
bounded pseudoconvex and fat
but which do not satisfy \thetag{$\ast$} (which gives the answer
about the existence of such domains in \cite{Jar-Pfl~2}).
This counterexample is found
in the class of Zalcman type domains, which were considered in the context
of \thetag{$\ast$} and Bergman completeness in \cite{Ohs~2} and \cite{Chen}.
It follows from the papers above that there are Zalcman type domains
satisfying \thetag{$\ast$} and being Bergman complete and non-hyperconvex.
We show
that there are Zalcman type domains, which neither satisfy \thetag{$\ast$}
nor are Bergman complete.

Let us fix a sequence (of pairwise different points)
$\{a_j\}_{j=1}^{\infty}\subset E$ and a closed
disc $B\subset E$ such that $a_j\not\in B$, $a_j\to 0$
and $0\in B$ (automatically $0\in\partial B$).

Below we shall consider only sequence of positive numbers $r_j$ such that
$\bar\triangle (a_j,r_j)\cap\bar\triangle(a_k,r_k)=\emptyset$
for any $j\neq k$ and $B\cap\bar\triangle(a_j,r_j)=\emptyset$.

\proclaim{Lemma 5.1}
We may choose $r_j$ so that there is a constant $M<\infty$
such that
$$
K_{D_N}(z)<M \text{ for any $z\in B$, $N=1,2,\ldots$},\tag{5.1}
$$
where $D_N:=E\setminus(\bigcup_{j=1}^N\bar\triangle(a_j,r_j))$.
\endproclaim
\demo{Proof} We define $r_j$ inductively. Since $E\setminus\bar
\triangle(a_1,r_1)$
increases to $E\setminus\{a_1\}$ as $r_1$ decreases to $0$
and $K_{E\setminus\{a_1\}}$ coincides with $K_E$ on $E\setminus\{a_1\}$
we get that there is a constant $M\in\Bbb R$ such that
$K_{E\setminus\bar\triangle(a_1,r_1)}<M$ on $B$ for sufficiently
small $r_1>0$.

Assume that we have already chosen $r_1,\ldots,r_N$ such that
$$
K_{D_N}<M \text{ on $B$}\tag{5.2}
$$
($D_N$ is defined as in the lemma).
Since $D_N\setminus\bar\triangle(a_{N+1},r_{N+1})$ increases to
$D_N\setminus\{a_{N+1}\}$ and $K_{D_N\setminus\{a_{N+1}\}}$ coincides with
$K_{D_N}$ on $D_N\setminus\{a_{N+1}\}$, we get as previously
(use \thetag{5.2})
that $K_{D_N\setminus\bar\triangle(a_{N+1},r_{N+1})}<M$ on $B$ for sufficiently
small $r_{N+1}>0$, which completes the proof.
\qed
\enddemo

\proclaim{Proposition 5.2} There is a sequence $s_j\to 0$, $0<s_j\leq r_j$
and a domain
$G:=E\setminus(\bigcup_{j=1}^{\infty}\bar\triangle(a_j,s_j)\cup \{0\})$
satisfying the property
$$
K_G(z)<M,\quad z\in B\cap G.
$$
\endproclaim
\demo{Proof} Let us fix an increasing sequence of compact sets $L_N$
such that $\bigcup_{N=1}^{\infty}L_N=\operatorname{int}B$.

We claim that we may choose a family of positive numbers
$\{s_N^j\}_{1\leq N<\infty,N\leq j}$ such that
$s^j_{N+1}\leq s^j_N\leq r_j$ for $j\geq N+1$ and $s^N_N\leq r_N$
such that for the domain
$$
G_N:=E\setminus
(\bigcup_{j=N+1}^{\infty}\bar\triangle(a_j,s_N^j)\cup\bigcup_{j=1}^N
\bar\triangle(a_j,s_j^j)\cup\{0\})
$$
we have $K_{G_N}<M$ on $L_N$.

Assume for a while that such a choice can be done. Then define $s_j:=s_j^j$.
Since $G_N\subset G$, we have $K_{G_N}\geq K_G$ for any $N$,
in particular, for any $N$ $K_G<M$ on $L_N$,
which completes the proof.

We define the desired family inductively with respect to $N$.
Let $0<t<1$.
Since
$E\setminus(\bigcup_{j=2}^{\infty}\bar\triangle(a_j,tr_j)\cup
\bar\triangle(a_1,r_1)\cup\{0\})$ increases to
$E\setminus(\bigcup_{j=2}^{\infty}\{a_j\}\cup\bar\triangle(a_1,r_1)\cup\{0\})$
as $t$ decreases to $0$
and the Bergman kernel of the last domain is the restriction
to this domain of $K_{D_1}$, we get that for $t$ sufficiently small
$K_{G_1}<M$ on $L_1$, where $s^1_1:=r_1$, $s^j_1=tr_j$, $j\geq 2$.

Assume that the construction has been succesful for $N$
(i.e. we have defined already all $s^j_k$,
$j\geq k$, $k\leq N$). Let $0<t<1$.
Since $E\setminus(\bigcup_{j=N+2}^{\infty}\bar\triangle(a_j,ts^j_N)\cup
\bar\triangle(a_{N+1},s^{N+1}_N)\cup\bigcup_{j=1}^N\bar\triangle(a_j,s_j^j)
\cup\{0\})$ increases to
$E\setminus(\bigcup_{j=N+2}^{\infty}\{a_j\}\cup
\bar\triangle(a_{N+1},s^{N+1}_N)\cup\bigcup_{j=1}^N\bar\triangle(a_j,s_j^j)
\cup\{0\})$ as $t$ decreases to $0$
and the Bergman kernel of the last domain is the restriction
to this domain of $K_{E\setminus(\bar\triangle(a_{N+1},s^{N+1}_N)
\cup\bigcup_{j=1}^N\bar\triangle(a_j,s^j_j))}$
(which is smaller than or equal to
$K_{D_{N+1}}$ because $D_{N+1}$ is a subset of the considered domain),
we get that defining for $t$ sufficiently small
$s^{N+1}_{N+1}:=s^{N+1}_N$,
$s^j_{N+1}:=ts^j_N$, $j\geq N+2$ the inequality
$K_{G_{N+1}}<M$ holds on $L_{N+1}$.
\qed
\enddemo

Let us remark that because of the property $K_{D_j}(z)\to K_D(z)$
locally uniformly
for any sequence $\{D_j\}_{j=1}^{\infty}$ of domains such that
$D_j\subset D_{j+1}$ and $\bigcup_{j=1}^{\infty}D_j=D$
($D$ is a bounded domain) we get easily that $\beta_{D_j}\to\beta_D$
locally uniformly on $D\times\Bbb C^n$ (although the convergence
in contrast to the convergence of Bergman kernels need not be monoton).

Based on the above property of the Bergman kernel we present below
a similar construction (to that from Proposition 5.2) leading
to a domain having the assumptions as in Proposition 5.2 and,
additionally, not Bergman complete. We denote
$\beta_{D}(z):=\beta_D(z;1)$.

\proclaim{Lemma 5.3} There are a constant $M_1\in\Bbb R$ and
a family of tuples
$\Lambda=\bigcup_{N=0}^{\infty}\Lambda_N$,
where:
$$
\Lambda_0=\emptyset,\;\Lambda_N\subset(0,s_1]\times\ldots\times(0,s_N];
$$
for any $N$ if $\lambda\in\Lambda_N$ then there is some
$\lambda_{N+1}$ such that
for any $0<s\leq\lambda_{N+1}$ $(\lambda,s)\in\Lambda_{N+1}$;
\newline
for any $\lambda=(\lambda_1,\ldots,\lambda_N)\in\Lambda$ we have that
$\beta_{D_{\lambda}}<M_1$ on $B$, where
$$
D_{\lambda}:=E\setminus\bigcup_{j=1}^N\bar\triangle(a_j,\lambda_j).
$$
\endproclaim
\demo{Proof} The proof goes similarily as that of Lemma 5.1.
We proceed using induction.
Since $E\setminus\bar\triangle(a_1,t)$
increases to $E\setminus\{a_1\}$ as $t$ decreases to $0$
and $\beta_{E\setminus\{a_1\}}$ coincides with $\beta_E$
on $E\setminus\{a_1\}$
we get that there is a constant $M_1\in\Bbb R$ such that
$\beta_{E\setminus\bar\triangle(a_1,t)}<M_1$ on $B$ for any
$0<t\leq \lambda_1\leq s_1$. We define $\Lambda_1:=(0,\lambda_1]$.

Assume that we have already defined $\Lambda_1,\ldots,\Lambda_N$ such that
Lemma is satisfied, in particular,
$$
\beta_{D_{\lambda}}<M_1 \text{ on $B$}
$$
for any $\lambda\in\Lambda_N$.

Fix any $\lambda\in\Lambda_N$.
Since $D_{\lambda}\setminus\bar\triangle(a_{N+1},t)$ increases to
$D_{\lambda}\setminus\{a_{N+1}\}$ and $\beta_{D_{\lambda}
\setminus\{a_{N+1}\}}$ coincides with
$\beta_{D_{\lambda}}$ on $D_{\lambda}\setminus\{a_{N+1}\}$,
we get as previously
that $\beta_{D_{\lambda}\setminus\bar\triangle(a_{N+1},t)}<M_1$
on $B$ for sufficiently small $t>0$, which completes the proof.
\qed
\enddemo

\proclaim{Proposition 5.4} There is a sequence $\lambda_j\to 0$,
$0<\lambda_j\leq s_j$ and a domain
$G:=E\setminus(\bigcup_{j=1}^{\infty}\bar\triangle(a_j,\lambda_j)\cup \{0\})$
satisfying the property
$$
\beta_G(z)\leq M_1,\quad z\in B\cap G.
$$
\endproclaim
\demo{Proof} Let us fix an increasing sequence of compact sets $L_N$
such that $\bigcup_{N=1}^{\infty}L_N=\operatorname{int}B$.

Without loss of generality we may assume that $s_1=\lambda^1$
($s_1$ is from Proposition 5.2 and $\lambda^1$ from Lemma 5.3).

It is sufficient to find sequences $\{\lambda^j\}_{j=1}^{\infty}$
and $\{t_j\}_{j=1}^{\infty}\subset(0,1)^{\Bbb N}$
such that $\lambda^N\in\Lambda_N$, $\lambda_N=t_1\cdot\ldots\cdot t_Ns_N$,
$\lambda^{N+1}=(\lambda^N,\lambda_{N+1})$,
and $\beta_{D_N}<M_1$ on $L_N$, where $D_N:=D_{\lambda^N}\setminus
\bigcup_{j=N+1}^{\infty}\bar\triangle(a_j,t_1\cdot\ldots\cdot t_Ns_j)$.

Put $\lambda_1:=\lambda^1$($=s_1$), $t_1:=1$.
Then for $1>t>0$ small enough the Bergman
metric on
$D_{\lambda^1}\setminus(\bigcup_{j=2}^{\infty}\bar\triangle(a_j,ts_j)
\cup\{0\})$
is less than $M_1$ on $L^1$ for $0<t\leq t_2<1$, we may also assume that
$\lambda^2:=(\lambda^1,t_2s_2)\in\Lambda^2$.

Assume that the construction has been succesful for $N$
(i.e. we have defined already all $t_j$, $j=1,\ldots,N$
and $\lambda^j$, $j=1,\ldots,N$).
Let $0<t<1$.
Since $D_{\lambda^N}\setminus(\bigcup_{j=N+1}^{\infty}\bar\triangle
(a_j,tt_1\cdot\ldots\cdot t_Ns_j)
\cup\{0\})$ increases to
$D_{\lambda^N}\setminus(\bigcup_{j=N+1}^{\infty}\{a_j\}
\cup\{0\})$ as $t$ decreases to $0$
and the Bergman metric of the last domain is the restriction
to this domain of $\beta_{D_{\lambda^N}}$
we may choose $t_{N+1}$ and then define $\lambda^{N+1}:=(\lambda^N,
t_1\cdot\ldots\cdot t_{N+1}s_{N+1})\in\Lambda^{N+1}$
having the desired properties.
\qed
\enddemo

\subheading{Remark 5.5} Note that the above mentioned results
may be put in some more general context. Two principal
properties that were used were the following both: $K_D$ and
$\beta_D$ do not change after deleting a discrete subset and both are
continuous with respect to the increasing family of domains. Applying
the same procedure we may prove for instance that there are Zalcman type
domains, which are not Carath\'eodory complete. Consequently, there are
Zalcman type domains,
without peak functions in $0$
(see \cite{Rud}).

\enddocument